\documentclass[psamsfonts,reqno]{amsart}
\usepackage{amscd,amssymb,eucal,latexsym,xy}
\xyoption{all}
\xyoption{ps}

\newcommand{\jj}{\vee}
\newcommand{\mm}{\wedge}
\newcommand{\JJ}{\bigvee}

\newcommand{\uu}{\cup}

\newcommand{\ci}{\subseteq}
\newcommand{\set}[1]{\{#1\}}
\newcommand{\setm}[2]{\{\,#1\mid#2\,\}}
\def\vv<#1>{\langle#1\rangle}

\newcommand{\ga}{\alpha}
\newcommand{\gb}{\beta}

\newcommand{\gd}{\delta}
\renewcommand{\ge}{\varepsilon}
\newcommand{\gf}{\varphi}

\newcommand{\gk}{\kappa}
\newcommand{\gl}{\lambda}

\newcommand{\gn}{\nu}
\newcommand{\go}{\omega}
\newcommand{\gp}{\pi}

\newcommand{\gr}{\varrho}

\newcommand{\gx}{\xi}
\newcommand{\gy}{\psi}

\newcommand{\tbf}{\textbf}
\newcommand{\tup}{\textup}

\DeclareMathAlphabet{\Bi}{OT1}{cmm}{b}{it}

\newcommand{\E}[1]{\mathcal{#1}}

\newcommand{\ol}[1]{\overline{#1}}

\def\con#1=#2(#3){#1\equiv#2\pod{#3}}
\providecommand{\bysame}{\makebox[3em]{\hrulefill}\thinspace}
\newcommand{\q}{\quad}
\newcommand{\qq}{\qquad}
\newcommand{\iso}{\cong}

\theoremstyle{plain}

\newtheorem*{mthm}{Theorem}
\newtheorem{lemma}{Lemma}[section]
\newtheorem{theorem}[lemma]{Theorem}
\newtheorem{proposition}[lemma]{Proposition}
\newtheorem{corollary}[lemma]{Corollary}

\newtheorem*{stat}{\name}
\newcommand{\name}{testing}

\theoremstyle{definition}
\newtheorem{definition}[lemma]{Definition}
\newtheorem{example}[lemma]{Example}
\newtheorem{problem}{Problem}

\theoremstyle{remark}

\newtheorem*{notation}{Notation}

\newenvironment{all}[1]{\renewcommand{\name}{#1}\begin{stat}}
                        {\end{stat}}

\newcommand{\qedc}{{\qed}~{\rm Claim~{\theclaim}.}}

\newcommand{\two}{\mathbf{2}}
\newcommand{\id}{\mathrm{id}}
\newcommand{\jz}{$\set{\jj,0}$}
\newcommand{\jzu}{$\set{\jj,0,1}$}

\newcommand{\Endk}{\mathrm{End}_K}
\newcommand{\NN}[2]{\mathrm{N}_{#1}({#2})}
\newcommand{\EE}[2]{\mathrm{E}_{#1}({#2})}

\newcommand{\famm}[2]{\langle\,#1\mid#2\,\rangle}
\newcommand{\sstrut}{\vrule height 5pt depth 1pt width 0pt}

\DeclareMathOperator{\Con}{Con}
\DeclareMathOperator{\Conc}{Con_c}
\DeclareMathOperator{\Id}{Id}
\DeclareMathOperator{\Idc}{Id_c}
\DeclareMathOperator{\im}{im}
\DeclareMathOperator{\M}{M}

 \author[F.~Wehrung]{Friedrich Wehrung}
 \address{CNRS, ESA 6081\\
          Universit\'e de Caen, Campus II\\
          D\'epartement de Math\'ematiques\\
          B.P. 5186\\
          14032 CAEN Cedex\\
          FRANCE}
 \email{wehrung@math.unicaen.fr}
 \urladdr{http://www.math.unicaen.fr/\~{}wehrung}


\begin{document}

\title[Representation of algebraic distributive lattices]
{Representation of algebraic distributive lattices 
with $\aleph_1$ compact elements as ideal
lattices of regular rings}

\begin{abstract}
We prove the following result:

\begin{mthm}
Every algebraic distributive lattice $D$ with at most
$\aleph_1$ compact elements is isomorphic to the ideal lattice of a
von~Neumann regular ring $R$.
\end{mthm}

(By earlier results of the author, the $\aleph_1$ bound is \emph{optimal}.)
Therefore, $D$ is also
isomorphic to the congruence lattice of a sectionally complemented modular
lattice $L$, namely, the principal right ideal lattice of $R$. Furthermore,
if the largest element of $D$ is compact, then one can assume that $R$ is
unital, respectively, that $L$ has a largest element. This extends
several known results of G.M. Bergman, A.P. Huhn, J. T\r uma,
and of a joint work of G.~Gr\"atzer, H.~Lakser, and the author, and
it solves Problem~2 of the survey paper \cite{GoWe1}.

The main tool used in the proof of our result is an amalgamation
theorem for semilattices and \emph{algebras} (over a given division ring), a
variant of previously known amalgamation theorems for semilattices and
\emph{lattices}, due to J. T\r uma, and G. Gr\"atzer, H. Lakser, and the
author.
\end{abstract}

\maketitle

\section*{Introduction}

It is a well-known and easy fact that the lattice of ideals of any
(von~Neumann) regular ring is algebraic and distributive.
In unpublished notes from 1986, G.M. Bergman \cite{Berg86} proves the
following converse of this result:

\begin{all}{Bergman's Theorem}
Every algebraic distributive lattice $D$ with countably many compact
elements is isomorphic to the ideal lattice of a regular
ring $R$, such that if the largest element of $D$ is compact, then $R$
is unital.
\end{all}

On the negative side, the author of the present paper,
in \cite{Weurp}, using his construction in
\cite{Wehr1}, proved that Bergman's Theorem cannot be extended to algebraic
distributive lattices with $\aleph_2$ many compact elements (or more). This
left a gap at the size $\aleph_1$, expressed by the statement of the
following problem:

\begin{all}{Problem~2 of \cite{GoWe1}}
Let $D$ be an algebraic distributive lattice with at most $\aleph_1$ compact
elements. Does there exist a regular ring $R$ such that the ideal lattice of
$R$ is isomorphic to $D$?
\end{all}

In this paper, we provide a positive solution to this problem, see
Theorem~\ref{T:Repr}. Of independent interest is an amalgamation result of
ring-theoretical nature, mostly inspired by the lattice-theoretical
constructions in \cite{TuSR} and \cite{GLWe}, see Theorem~\ref{T:Amalg}.
This result is the main tool used in the proof of Theorem~\ref{T:Repr}.

Once the Amalgamation Theorem (Theorem~\ref{T:Amalg}) is proved, the
representation result follows from standard techniques, based on the
existence of lattices called
\emph{$2$-frames} in \cite{Dobb86}, or \emph{lower finite
$2$-lattices} in \cite{Dito84}. Such a technique has, for
example, been used successfully in \cite{Huhn89a,Huhn89b}, where A.P. Huhn
proves that every distributive algebraic lattice $D$ with at most
$\aleph_1$ compact elements is isomorphic to the congruence lattice of a
lattice $L$. Theorem~\ref{T:Repr} provides a strengthening of Huhn's
result, namely, it makes it possible to have $L$ \emph{sectionally
complemented} and \emph{modular}, see Corollary~\ref{C:Repr}. Note that we
already obtained $L$ relatively complemented with zero (thus sectionally
complemented), though not modular, in \cite{GLWe}.

We do not claim any originality about the proof methods used in this paper.
Most of what we do amounts to translations between known concepts and
proofs in universal algebra, lattice theory, and ring theory. However, the
interconnections between these domains, as they are, for example, presented
in \cite{GoWe1}, are probably not well-established enough to trivialize the
results of this paper.

\section{Basic concepts}\label{S:basic}

\subsection*{Lattices, semilattices}

References for this section are \cite{Birk67,Grat71,GLT2}.

Let $L$ be a lattice. We say that $L$ is \emph{complete}, if every subset of
$L$ has a supremum.
An element $a$ of $L$ is \emph{compact}, if for every subset $X$
of $L$ such that the supremum of $X$, $\JJ X$, exists, $a\leq\JJ X$ implies
that there exists a finite subset $Y$ of $X$ such that $a\leq\JJ Y$.
The \emph{unit} of a lattice is its largest element, if it exists.

We say that $L$ is \emph{algebraic}, if it is complete and every element is
the supremum of compact elements. If $L$ is an algebraic lattice, then the
set $S$ of all compact elements of $L$ is closed under the join operation,
and it contains $0$ (the least element of $L$) as an element. We say that
$S$ is a \emph{\jz-semilattice}, that is, a commutative, idempotent monoid
(the monoid operation is the join).

If $S$ is a \jz-semilattice, an \emph{ideal} of $S$ is a nonempty,
hereditary subset of $S$ closed under the join operation.
The set $\Id S$ of all ideals of $S$,
partially ordered by containment, is an algebraic lattice, and the
compact elements of $\Id S$ are exactly the principal ideals
$(s]=\setm{t\in S}{t\leq s}$, for $s\in S$. In particular, the
semilattice of all compact elements of $\Id S$ is isomorphic to $S$.
Conversely, if $L$ is an algebraic lattice and if $S$ is the semilattice of
all compact elements of $L$, then the map from $L$ to $\Id S$ that with every
element $x$ of $L$ associates $\setm{s\in S}{s\leq x}$ is an isomorphism.
It follows that algebraic lattices and \jz-semilattices are categorically
equivalent. The class of homomorphisms of algebraic lattices that
correspond, through this equivalence, to \jz-homomorphisms of semilattices,
are the \emph{compactness preserving, $\JJ$-complete homomorphisms} of
algebraic lattices. (We say that a homomorphism $f\colon A\to B$ of algebraic
lattices is \emph{$\JJ$-complete}, if $\JJ f[X]=f\left(\JJ X\right)$, for
every (possibly empty) subset $X$ of $A$.) We observe that if $B$ is finite,
then any homomorphism from $A$ to $B$ is compactness-preserving.

A \jz-semilattice $S$ is \emph{distributive}, if its ideal lattice $\Id S$
is a distributive lattice. Equivalently, $S$ satisfies the following
statement:
 \[
 (\forall a,b,c)\bigl(c\leq a\jj b\Rightarrow(\exists x,y)
 (x\leq a\text{ and }y\leq b\text{ and }c=x\jj y)\bigr).
 \]
For a lattice $L$, the set $\Con L$ of all congruences of
$L$, endowed with containment, is an algebraic lattice. Its semilattice of
compact elements is traditionally denoted by $\Conc L$, the semilattice of
\emph{finitely generated congruences} of $L$.

We say that $L$ is
\begin{itemize}
\item[---] \emph{modular}, if $x\mm(y\jj z)=(x\mm y)\jj z$ for all $x$, $y$,
$z\in L$ such that $x\geq z$;

\item[---] \emph{complemented}, if it has a least element, denoted by $0$,
a largest element, denoted by $1$, and for all $a\in L$, there exists
$x\in L$ such that $a\mm x=0$ and $a\jj x=1$;

\item[---] \emph{sectionally complemented}, if it has a least element,
denoted by $0$, and for all $a$, $b\in L$ such that $a\leq b$, there exists
$x\in L$ such that $a\mm x=0$ and $a\jj x=b$.
\end{itemize}

We denote by $\two$ the two-element lattice.

\subsection*{Rings, algebras}

All the rings encountered in this work are associative, but not necessarily
unital. A ring $R$ is \emph{regular} (in von~Neumann's sense), if it
satisfies the statement $(\forall x)(\exists y)(xyx=x)$. If $R$ is a
regular ring, then the set of all principal right ideals of $R$, partially
ordered by containment, is a \emph{sectionally complemented modular
lattice}, see, for example, \cite[Page 209]{FrHa56}.

If $R$ is a ring, then we denote by $\Id R$ the set of all two-sided ideals
of $R$, partially ordered by containment. Then
$\Id R$ is an algebraic modular lattice, which turns out to be distributive
if $R$ is regular. We denote by $\Idc R$ the semilattice of all compact
elements of $\Id R$, that is, the finitely generated two-sided ideals of
$R$. It is to be noted that $\Idc$ can be extended to a \emph{functor} from
rings and ring homomorphisms to \jz-semilattices and \jz-homomorphisms, and
that this functor preserves direct limits.

If $K$ is a division ring, a \emph{$K$-algebra} is a ring $R$ endowed with
a structure of two-sided vector space over $K$ such that the equalities
 \[
 \gl(xy)=(\gl x)y,\quad (x\gl)y=x(\gl y),\quad (xy)\gl=x(y\gl)
 \]
hold for all $x$, $y\in R$, and $\gl\in K$. Such a structure is called a
\emph{$K$-ring} in \cite[Section~1]{Cohn59}. Most of the rings that we shall
encounter in this work are, in fact, algebras.

\section{Embedding into V-simple algebras}

\begin{definition}\label{D:Vsimp}
A unital, regular ring $R$ is \emph{V-simple}, if $R_R$ is isomorphic to all
its nonzero principal right ideals, and there are nonzero principal right
ideals $I$ and $J$ of $R$ such that $I\oplus J=R_R$.
\end{definition}

It is obvious that if $R$ is V-simple, then it is simple. The converse is
obviously false, for example, if $R$ is a field.

\begin{notation}
Let $\gk$ be an infinite cardinal number,
let $U$ be a two-sided vector space over a division ring $K$.
For example, if $I$ is any set, then the set $K^{(I)}$ of all
$I$-families with finite support of elements of $K$ is endowed with
a natural structure of two-sided vector space over $K$.

Let $\NN{\gk}{U}$ be the subset of the algebra $\Endk(U)$
of \emph{right} $K$-vector space endomorphisms of $U$ defined by
 \[
 \NN{\gk}{U}=\setm{f\in\Endk(U)}{\dim_K\im f<\gk}.
 \]
It is obvious that $\NN{\gk}{U}$ is a two-sided ideal of $\Endk(U)$.
We define a $K$-algebra, $\EE{\gk}{U}$, by
 \[
 \EE{\gk}{U}=\Endk(U)/\NN{\gk}{U}.
 \]
\end{notation}

The idea behind the proof of Lemma~\ref{L:Endgk} and
Proposition~\ref{P:EmbVsimple} is old, see, for example, Theorem~3.4 in
\cite{Ara}. For convenience, we recall the proofs here.

\begin{lemma}\label{L:Endgk}
Let $\gk$ be an infinite cardinal number, let $U$ be a two-sided vector
space of dimension $\gk$ over a division ring $K$. Then $\EE{\gk}{U}$ is
a unital, regular, V-simple $K$-algebra.
\end{lemma}

\begin{proof}
Note that $\EE{\gk}{U}$ is nontrivial, because $\gk\leq\dim_KU$.
Since the endomorphism ring $\Endk(U)$ is regular, so is also the quotient
ring $\EE{\gk}{U}=\Endk(U)/\NN{\gk}{U}$, see Lemma~1.3 in
\cite{GvnRR}.

Furthermore, the principal right ideals of $R=\Endk(U)$ are exactly the
ideals of the form
 \[
 I_X=\setm{f\in R}{\im f\ci X},
 \]
for a subspace $X$ of $U$. If $X$ and $Y$ are subspaces of $U$, then
$X\iso Y$ implies that $I_X\iso I_Y$.
If $X$ is a subspace of dimension $\gk$ of $U$, then
$X$ can be decomposed as $X=X_0\oplus X_1$, where
$\dim_KX_0=\dim_KX_1=\gk$, hence
 \[
 [I_X]=[I_{X_0}]+[I_{X_1}]=2[I_X],
 \]
where $[I]$ denotes the isomorphism class of a right ideal $I$.
Since $X\iso U$, $[I_X]=[I_U]=[R]$.
However, if $X$ is a subspace of $U$ of dimension $<\gk$, then the image of
$I_X$ in $\EE{\gk}{U}$ is the zero ideal. The conclusion follows.
\end{proof}

\begin{proposition}\label{P:EmbVsimple}
Let $K$ be a division ring.
Every unital $K$-algebra has a unital embedding into a
unital, regular, V-simple $K$-algebra.
\end{proposition}

\begin{proof}
Let $R$ be a unital $K$-algebra. Put $\gk=\aleph_0+\dim_KR$,
where $\dim_KR$ denotes the \emph{right dimension} of $R$ over $K$,
and $U=R^{(\gk)}$, the $R$-algebra of all $\gk$-sequences with finite
support of elements of $R$. We put
 \[
 S=\EE{\gk}{U}.
 \]
Since the dimension of $U$ over $K$ equals $\gk$, it follows from
Lemma~\ref{L:Endgk} that $S$ is a unital, regular, V-simple $K$-algebra.

Define a map $\gf\colon R\to\Endk(U)$, by the rule
 \[
 \gf(a)\colon U\to U,\quad x\mapsto ax,
 \]
for all $a\in R$. It is easy to see that $\gf$ is a unital ring homomorphism
from $R$ to $\Endk(U)$.

Let $a\in R\setminus\set{0}$.
If $\famm{e_\gx}{\gx<\gk}$ denotes the canonical basis of $U$ over
$R$, then the range of $\gf(a)$ contains all the elements $ae_\gx$, for
$\gx<\gk$, thus its dimension over $K$ is greater than or equal to
$\gk$. Hence,
$\gf(a)$ does not belong to $\NN{\gk}{U}$. Therefore, the map $\gy$ from
$R$ to $S$ defined by the rule
 \[
 \gy(a)=\gf(a)+\NN{\gk}{U},
 \]
for all $a\in R$, is a unital $K$-algebra embedding from $R$ into $S$.
\end{proof}

\section{Amalgamation of algebras over a division ring}

The following fundamental result has been proved by P.M. Cohn, see
Theorem~4.7 in \cite{Cohn59}.
We also refer the reader to the outline presented in \cite[Page 110]{KMPT},
in the section ``Regular rings: AP''.

\begin{theorem}\label{T:amalgalg}
Let $K$ be a division ring, let $R$, $A$, and $B$ be unital $K$-algebras,
with $R$ regular.
Let $\ga\colon R\hookrightarrow A$ and $\gb\colon R\hookrightarrow B$ be
unital embeddings. Then there exist a
unital, regular $K$-algebra $C$,
$\ga'\colon A\hookrightarrow C$, and $\gb'\colon B\hookrightarrow C$,
such that $\ga'\circ\ga=\gb'\circ\gb$.
\end{theorem}

By combining Theorem~\ref{T:amalgalg} with Proposition~\ref{P:EmbVsimple},
we obtain immediately the following slight strengthening of
Theorem~\ref{T:amalgalg}:

\begin{lemma}\label{L:amalgalg2}
Let $K$ be a division ring, let $R$, $A$, and $B$ be unital $K$-algebras,
with $R$ regular.
Let $\ga\colon R\hookrightarrow A$ and $\gb\colon R\hookrightarrow B$ be
unital embeddings. Then there exist a unital, regular, V-simple $K$-algebra
$C$, unital embeddings $\ga'\colon A\hookrightarrow C$, and
$\gb'\colon B\hookrightarrow C$, such that $\ga'\circ\ga=\gb'\circ\gb$.
\end{lemma}

The following example partly illustrates the underlying complexity of
Theorem~\ref{T:amalgalg}, by showing that even in the case where $A$ and $B$
are finite-dimensional over $R$, one may \emph{not} be able to find a
finite-dimensional solution $C$ to the amalgamation problem:

\begin{example}\label{Eq:NonFin}
Let $K$ be any division ring. We construct unital, matricial extensions $A$
and $B$ of the regular $K$-algebra $R=K^2$ such that the amalgamation
problem of $A$ and $B$ over $R$ has no finite-dimensional solution.
\end{example}

\begin{proof}
We put $A=\M_2(K)$ (resp., $B=\M_3(K)$), the ring of all square matrices of
order $2$ (resp., $3$) over $K$, endowed with their canonical
$K$-algebra structures. We define unital embeddings of
$K$-algebras $f\colon R\hookrightarrow A$ and
$g\colon R\hookrightarrow B$ as follows:
 \begin{equation}\label{Eq:Matr}
 f(\vv<x,y>)=
  \begin{pmatrix}
  x&0\\
  0&y
  \end{pmatrix},
 \quad\text{and}\quad
 g(\vv<x,y>)=
  \begin{pmatrix}
  x&0&0\\
  0&x&0\\
  0&0&y
  \end{pmatrix},
 \end{equation}
for all $x$, $y\in K$.

Suppose that the amalgamation problem of $A$ and $B$ over $R$ (with respect
to $f$ and $g$) has a solution, say, $f'\colon A\hookrightarrow C$ and
$g'\colon B\hookrightarrow C$, where $C$ is a \emph{finite-dimensional}, not
necessarily unital $K$-algebra and $f'$, $g'$ are embeddings of unital
$K$-algebras. For a positive integer $d$, we denote by $e_{i,j}^d$,
for $1\leq i,j\leq d$, the canonical matrix units of the matrix ring
$\M_d(K)$. So,
\eqref{Eq:Matr} can be rewritten as
 \begin{equation}\label{Eq:Matr2}
 f(\vv<x,y>)=e_{1,1}^2x+e_{2,2}^2y,
 \quad\text{and}\quad
 g(\vv<x,y>)=(e_{1,1}^3+e_{2,2}^3)x+e_{3,3}^3y,
 \end{equation}
for all $x$, $y\in K$. Put $u_{i,j}=f'(e_{i,j}^2)$, for all $i$,
$j\in\set{1,2}$, and
$v_{i,j}=g'(e_{i,j}^3)$, for all $i$, $j\in\set{1,2,3}$. Then apply $f'$
(resp.
$g'$) to the first (resp., second) equality of \eqref{Eq:Matr2}. Since
$f'\circ f=g'\circ g$, we obtain that the equality
 \[
 u_{1,1}x+u_{2,2}y=(v_{1,1}+v_{2,2})x+v_{3,3}y,
 \]
holds, for all $x$, $y\in K$. Specializing to $x$, $y\in\set{0,1}$ yields
the equalities
 \begin{align}
 u_{1,1}&=v_{1,1}+v_{2,2},\label{Eq:uv12}\\
 u_{2,2}&=v_{3,3}.\label{Eq:u0v0}
 \end{align}
However, the elements $u_{i,j}$ of $C$, for $i$, $j\in\set{1,2}$, satisfy
part of the equalities defining matrix units in $C$, namely,
$u_{i,j}u_{k,l}=\gd_{j,k}u_{i,l}$, for all $i$, $j$, $k$,
$l\in\set{1,2}$ ($\gd_{-}$ denotes here the Kronecker symbol).
Similarly, the elements $v_{i,j}$ of~$C$, for
$i$, $j\in\set{1,2,3}$, satisfy the equalities
$v_{i,j}v_{k,l}=\gd_{j,k}v_{i,l}$, for all $i$, $j$, $k$,
$l\in\set{1,2,3}$. Therefore, by \eqref{Eq:uv12} and \eqref{Eq:u0v0},
$\dim_K(u_{1,1}C)=\dim_K(v_{1,1}C)+\dim_K(v_{2,2}C)=2\dim_K(v_{3,3}C)=
2\dim_K(u_{2,2}C)=2\dim_K(u_{1,1}C)$. But then, since $C$ is
finite-dimensional, $\dim_K(u_{1,1}C)=0$, so
$\dim_K(f'A)=\dim_K(g'B)=0$, a contradiction.
\end{proof}

For a further discussion of Example~\ref{Eq:NonFin}, see the comments
following the statement of Problem~\ref{Pb:LocMat} in Section~\ref{S:Pbs}.

\section{The amalgamation theorem}

\begin{definition}\label{D:VBool}
Let $K$ be a division ring. A $K$-algebra $R$ is \emph{V-Boolean}, if it is
isomorphic to a finite direct product of V-simple $K$-algebras.
\end{definition}

\begin{theorem}\label{T:Amalg}
Let $K$ be a division ring, let $R_0$, $R_1$, and $R_2$ be
unital $K$-algebras, with $R_0$ regular, let $S$ be a finite Boolean
lattice. For $k\in\set{1,2}$, let
$f_k\colon R_0\to R_k$ be a homomorphism of unital
$K$-algebras and let $\gy_k\colon\Id R_k\to S$ be a unit-preserving
$\JJ$-complete homomorphism, such that
$\gy_1\circ\Id f_1=\gy_2\circ\Id f_2$. Then there exist
a unital, regular, V-Boolean $K$-algebra $R$, homomorphisms of unital
$K$-algebras $g_k\colon R_k\to R$, for $k\in\set{1,2}$, and an isomorphism
$\ga\colon \Id R\to S$ such that $g_1\circ f_1=g_2\circ f_2$
and $\ga\circ\Id g_k=\gy_k$ for $k\in\set{1,2}$ (see Figure~\tup{1}).
\end{theorem}
 \[
{
\xymatrixrowsep{2pc}\xymatrixcolsep{1.5pc}
\def\labelstyle{\displaystyle}
\xymatrix{
& & & & \Id R\ar[d]|-{\sstrut\ga} &\\
& R & & & S \\
R_1\ar[ru]^{g_1} & & R_2\ar[lu]_{g_2}
& \Id R_1\ar[ur]|-{\gy_1}\ar[uur]|-{\Id g_1} & &
\Id R_2\ar[ul]|-{\gy_2}\ar[uul]|-{\Id g_2}\\
& R_0\ar[ul]^{f_1}\ar[ur]_{f_2} & &
\save+<0ex,-5ex>\drop{\text{Figure 1}}\restore &
\Id R_0\ar[ul]^{\Id f_1}\ar[uu]|-{\gy_0}\ar[ur]_{\Id f_2}
}
}
 \]
\begin{proof}
We adapt to the context of regular rings the lattice-theoretical
proof of Theorem~1 of \cite{GLWe}.

We put $\gy_0=\gy_1\circ\Id f_1=\gy_2\circ\Id f_2$. 
We start with the following case:
\smallskip

\noindent{\textbf{\textit{Case 1.}}} $S\iso\two$.
\smallskip

We put
$I_k=\setm{x\in R_k}{\gy_k(R_kxR_k)=0}$, for $k\in\set{0,1,2}$. Since
$\gy_k$ is a $\JJ$-complete homomorphism, $I_k$ is the largest ideal of
$R_k$ whose image under $\gy_k$ is zero. Furthermore, since $\gy_k$ is
unit-preserving, $I_k$ is a \emph{proper} ideal of $R_k$.

Next, we put $\ol{R}_k=R_k/I_k$, and we denote by
$p_k\colon R_k\twoheadrightarrow\ol{R}_k$ the canonical projection.

For $k\in\set{1,2}$, the equivalence $x\in I_0\Leftrightarrow f_k(x)\in I_k$
holds for all $x\in R_0$, thus there exists a unique unital \emph{embedding}
$\ol{f}_k\colon\ol{R}_0\hookrightarrow\ol{R}_k$ such that
$p_k\circ f_k=\ol{f}_k\circ p_0$, see Figure~2.
 \[
{
\def\labelstyle{\displaystyle}
\xymatrixrowsep{2pc}\xymatrixcolsep{1.5pc}
\xymatrix{
& & & R & \\
R_0\ar[r]^{f_k}\ar@{->>}[d]_{p_0} & R_k\ar@{->>}[d]^{p_k}
& \ol{R}_1\ar@{_(->}[ru]^{\ol{g}_1} &
& \ol{R}_2\ar@{^(->}[lu]_{\ol{g}_2}\\
\ol{R}_0\ar@{^(->}[r]_{\ol{f}_k} & \ol{R}_k &
\save+<0ex,-5ex>\drop{\text{Figure 2}}\restore
& \ol{R}_0\ar@{_(->}[ul]^{\ol{f}_1}\ar@{^(->}[ur]_{\ol{f}_2} &
}
}
 \]
Since $\ol{R}_0$ is regular and $\ol{f}_0$, $\ol{f}_1$ are unital
embeddings, there exist, by Lemma~\ref{L:amalgalg2}, a 
unital, regular, V-simple $K$-algebra $R$ and embeddings
$\ol{g}_k\colon\ol{R}_k\hookrightarrow R$, for $k\in\set{1,2}$, such that
$\ol{g}_1\circ\ol{f}_1=\ol{g}_2\circ\ol{f}_2$, see Figure~2.

We put $g_k=\ol{g}_k\circ p_k$, for $k\in\set{1,2}$, see Figure~3.
 \[
{
\def\labelstyle{\displaystyle}
\xymatrixrowsep{2pc}\xymatrixcolsep{1.5pc}
\xymatrix{
& R & \\
\ol{R}_1\ar@{_(->}[ru]^{\ol{g}_1} & & \ol{R}_2\ar@{^(->}[lu]_{\ol{g}_2} \\
R_1\ar@{->>}[u]|-{p_1}\ar@/^4pc/[uur]|-{g_1} &
\ol{R}_0\ar@{_(->}[ul]_{\ol{f}_1}\ar@{^(->}[ur]^{\ol{f}_2}
& R_2\ar@{->>}[u]|-{p_2}\ar@/_4pc/[uul]|-{g_2} \\
& R_0\ar@{->>}[u]|-{p_0}\ar[lu]^{f_1}\ar[ru]_{f_2}
\save+<0ex,-5ex>\drop{\text{Figure 3}}\restore & 
}
}
 \]
Therefore,
 \[
 g_1\circ f_1=\ol{g}_1\circ p_1\circ f_1=\ol{g}_1\circ\ol{f}_1\circ p_0
 =\ol{g}_2\circ\ol{f}_2\circ p_0=g_2\circ f_2.
 \]
Since $R$ is V-simple, it is simple, thus, since $S\iso\two$, there exists a
unique isomorphism $\ga\colon\Id R\to S$. To verify that $\ga\circ\Id
g_k=\gy_k$, for $k\in\set{1,2}$, it suffices to verify that $(\Id g_k)(I)=0$
if{f} $\gy_k(I)=0$, for every $I\in\Id R_k$. We proceed:
 \begin{align*}
 (\Id g_k)(I)=0&\q\text{if{f}}\q g_k[I]=0\\
 &\q\text{if{f}}\q\ol{g}_k\circ p_k[I]=0\\
 &\q\text{if{f}}\q p_k[I]=0\\
\intertext{(because $\ol{g}_k$ is an embedding)}
 &\q\text{if{f}}\q I\ci I_k\\
 &\q\text{if{f}}\q \gy_k(I)=0,
 \end{align*}
which concludes Case~1.

\smallskip
\noindent{\textbf{\textit{Case 2.}}} General case, $S$ finite Boolean.
\smallskip

Without loss of generality, $S=\two^n$, with $n<\go$. For $i<n$, let
$\gp_i\colon S\twoheadrightarrow\two$ be the projection on the $i$-th
coordinate. We apply Case~1 to the maps $\gp_i\gy_k$, for $k\in\set{1,2}$,
see Figure~4.
 \[
{
\xymatrixrowsep{2pc}\xymatrixcolsep{1.5pc}
\def\labelstyle{\displaystyle}
\xymatrix{
& & & & \Id R^{(i)}\ar[d]|-{\sstrut\ga_i} &\\
& R^{(i)} & & & \two \\
R_1\ar[ru]^{g_{1,i}} & & R_2\ar[lu]_{g_{2,i}}
& \Id R_1\ar[ur]|-{\pi_i\gy_1}\ar[uur]|-{\Id g_{1,i}} & &
\Id R_2\ar[ul]|-{\pi_i\gy_2}\ar[uul]|-{\Id g_{2,i}}\\
& R_0\ar[ul]^{f_1}\ar[ur]_{f_2} & &
\save+<0ex,-5ex>\drop{\text{Figure 4}}\restore &
\Id R_0\ar[ul]^{\Id f_1}\ar[uu]|-{\gy_0}\ar[ur]_{\Id f_2}
}
}
 \]
We obtain a unital, regular, V-simple $K$-algebra $R^{(i)}$, $K$-algebra
homomorphisms $g_{k,i}$, for $k\in\set{1,2}$, an isomorphism
$\ga_i\colon\Id R^{(i)}\to\two$, such that
$g_{1,i}\circ f_1=g_{2,i}\circ f_2$ and
$\ga_i\circ\Id g_{k,i}=\gp_i\circ\gy_k$, for $k\in\set{1,2}$,
see Figure~4.

Now we put $R=\bigoplus_{i<n}R^{(i)}$, with the componentwise ring
structure. So $R$ is a unital, regular, V-Boolean $K$-algebra. For
$k\in\set{1,2}$, we define a unital homomorphism $g_k\colon R_k\to R$
by the rule
 \[
 g_k(x)=\famm{g_{k,i}(x)}{i<n},\qq\text{for all }x\in R_k.
 \]
It is immediate that $g_1\circ f_1=g_2\circ f_2$.

Furthermore, observe that $\prod_{i<n}\Id R^{(i)}\iso\Id R$, \emph{via} the
isomorphism that sends a finite sequence $\famm{I_i}{i<n}$ to
$\bigoplus_{i<n}I_i$. Define an isomorphism $\ga\colon\Id R\to S$ by the rule
 \[
 \ga\left(\bigoplus_{i<n}I_i\right)=\famm{\ga_i(I_i)}{i<n},
 \qq\text{for all }I\in\Id R.
 \]
For $k\in\set{1,2}$ and $I\in\Id R_k$, we compute:
 \begin{align*}
 \ga\circ(\Id g_k)(I)&=\ga\left(\bigoplus_{i<n}(\Id g_{k,i})(I)\right)\\
 &=\famm{\ga_i\circ(\Id g_{k,i})(I)}{i<n}\\
 &=\famm{\gp_i\circ\gy_k(I)}{i<n}\\
 &=\gy_k(I),
 \end{align*}
so $\ga\circ\Id g_k=\gy_k$.
\end{proof}

\section{The representation theorem}\label{S:Repr}

We first state a useful lemma, see \cite{Dito84} and \cite{Dobb86}:

\begin{lemma}\label{L:ladder}
There exists a lattice $F$ of cardinality $\aleph_1$ satisfying the
following properties:
\begin{enumerate}
\item $F$ is \emph{lower finite}, that is, for all $a\in F$, the principal
ideal
 \[
 F[a]=\setm{x\in F}{x\leq a}
 \]
is finite.

\item Every element of $F$ has at most two immediate predecessors.
\end{enumerate}

\end{lemma}

\begin{theorem}\label{T:Repr}
Let $D$ be an algebraic distributive lattice with at most $\aleph_1$ compact
elements, and let $K$ be a division ring. Then there exists a
regular $K$-algebra $R$ satisfying the following properties:

\begin{enumerate}
\item $\Id R\iso D$.

\item If the largest element of $D$ is compact, then $R$ is a direct limit
of unital, regular, V-Boolean $K$-algebras and unital embeddings of
$K$-algebras. In particular, $R$ is unital.
\end{enumerate}

\end{theorem}

\begin{proof}
A similar argument has already been used, in different contexts,
in such various references as \cite{Dobb86,Huhn89a,Huhn89b,GLWe}.

We first translate the problem into the language of semilattices. This
amounts, by defining $S$ as the \jz-semilattice of compact elements of $D$,
to verifying the existence of a regular $K$-algebra $R$ satisfying the
following condition

\smallskip
(i') $\Idc R\iso S$,
\smallskip

\noindent
along with (ii). We do this first in the case where the largest element of
$D$ is compact, that is, where $S$ has a largest element. It is proved in
\cite{GoWe1} that $S$ is a direct limit of
finite Boolean \jz-semilattices and \jzu-homomorphisms, say,
 \[
 \vv<S,\gf_i>_{i\in I}=\varinjlim\vv<S_i,\gf^i_j>_{i\leq j\text{ in }I},
 \]
where $I$ is a directed partially ordered set, and
$\vv<S_i,\gf^i_j>_{i\leq j\text{ in }I}$ is a direct system of finite Boolean
\jz-semilattices and \jzu-homomorphisms (in particular,
$\gf^i_j\colon S_i\to S_j$ and $\gf_i\colon S_i\to S$,
for $i\leq j$ in $I$). Furthermore, one can take
$I$ countably infinite if $S$ is finite, and $|I|=|S|$ if $S$ is infinite.
In particular, $|I|\leq\aleph_1$.

Let $F$ be a lattice satisfying the conditions of Lemma~\ref{L:ladder}.
Since $|F|\geq|I|>0$, there exists a surjective map
$\gn_0\colon F\twoheadrightarrow I$. Since $F$ is lower finite, it is
well-founded, so we can define inductively an order-preserving, cofinal map
$\gn\colon F\to I$, by putting
 \[
 \gn(x)=\text{any element }i\text{ of }I\text{ such that }\gn_0(x)\leq i
 \text{ and }\gn(y)\leq i,\text{ for all }y<x,
 \]
for all $x\in F$. This is justified because $F$ is lower finite, and this
does not use part (ii) of Lemma~\ref{L:ladder}. As a conclusion, we see
that we may index our direct system by $F$ itself, that is,
\emph{we may assume that $I$ satisfies the conditions \tup{(i)},
\tup{(ii)} of Lemma~\tup{\ref{L:ladder}}}.

We shall now define inductively unital, regular, V-Boolean $K$-algebras
$R_i$, unital homomorphisms of $K$-algebras $f^i_j\colon R_i\to R_j$, and
isomorphisms $\ge_i\colon\Idc R_i\to S_i$, for $i\leq j$ in~$I$.

Let $\gr\colon I\to\go$ be the natural rank function, that is,
 \[
 \gr(i)=\sup\setm{\gr(j)}{j<i}+1,
 \]
for all $i\in I$.
For all $n<\go$, we put
 \[
 I_n=\setm{i\in I}{\gr(i)\leq n}.
 \]
By induction on $n<\go$, we construct V-Boolean $K$-algebras $R_i$
(note then that $\Id R_i=\Idc R_i$), maps
$\ge_i\colon\Id R_i\to S_i$, and unital homomorphisms of $K$-algebras
$f^i_j\colon R_i\to R_j$, for all $i\leq j$ in $I_n$, satisfying the
following properties:
\begin{itemize}
\item[(a)] $f^i_i=\id_{R_i}$, for all $i\in I_n$.

\item[(b)] $f^i_k=f^j_k\circ f^i_j$, for all $i\leq j\leq k$ in $I_n$.

\item[(c)] $\ge_i$ is a lattice isomorphism from $\Id R_i$ onto $S_i$,
for all $i\in I_n$.

\item[(d)] The following diagram is commutative, for all $i\leq j$ in $I_n$:
 \[
  \begin{CD}
  \Id R_i @>\Id f^i_j>> \Id R_j\\
  @V{\ge_i}VV @VV{\ge_j}V\\
  S_i @>>\gf^i_j> S_j
  \end{CD}
 \]
\end{itemize}

For $n=0$, it suffices to construct a V-Boolean $K$-algebra $R_0$ such that
$\Id R_0\iso S_0$. This is easy: if $p$ is the number of atoms of $S_0$,
take any V-simple $K$-algebra $R$, and put $R_0=R^p$.

Suppose having done the construction on $I_n$, we show how to extend it to
$I_{n+1}$. Let $i\in I_{n+1}$ such that $\gr(i)=n+1$. Denote by $i_0$ and
$i_1$ the two immediate predecessors of $i$ in $I$. Note that
$i_0$ and $i_1$ do not need to be distinct.
 \[
{
\xymatrixrowsep{2pc}\xymatrixcolsep{1.5pc}
\def\labelstyle{\displaystyle}
\xymatrix{
& & & & \Id R_i\ar[d]|-{\sstrut\ge_i} &\\
& R_i & & & S_i \\
R_{i_0}\ar[ru]^{g_0} & & R_{i_1}\ar[lu]_{g_1}
& \Id R_{i_0}\ar[ur]|-{\gf_i^{i_0}\ge_{i_0}}\ar[uur]|-{\Id g_0} & &
\Id R_{i_1}\ar[ul]|-{\gf_i^{i_1}\ge_{i_1}}\ar[uul]|-{\Id g_1}\\
& R_{i_0\mm i_1}\ar[ul]^{f_{i_0}^{i_0\mm i_1}}\ar[ur]_{f_{i_1}^{i_0\mm i_1}}
 & & \save+<0ex,-5ex>\drop{\text{Figure 5}}\restore &
\Id R_{i_0\mm i_1}\ar[ul]^{\Id f_{i_0}^{i_0\mm i_1}}
\ar[ur]_{\Id f_{i_1}^{i_0\mm i_1}}
}
}
 \]
By Theorem~\ref{T:Amalg}, there exist a unital, regular, V-Boolean $K$-algebra
$R_i$, unital homomorphisms of $K$-algebras $g_k\colon R_{i_k}\to R_i$, for
$k<2$, and an isomorphism $\ge_i\colon \Id R_i\to S_i$ such that the
following equalities hold:
 \begin{align}
 g_0\circ f^{i_0\mm i_1}_{i_0}&=g_1\circ f^{i_0\mm i_1}_{i_1};
 \label{Eq:g0fi0i1}\\
 \ge_i\circ\Id g_k&=\gf^{i_k}_i\ge_{i_k},\qq\text{for }k\in\set{1,2},
 \label{Eq:geiIdgk}
 \end{align}
see Figure~5. If $i_0=i_1$, we may replace $g_1$ by $g_0$: the diagrams
of Figure~5 remain commutative and \eqref{Eq:g0fi0i1}, \eqref{Eq:geiIdgk}
remain valid. Thus we may define $f^{i_0}_i=g_0$ and $f^{i_1}_i=g_1$, and
\eqref{Eq:g0fi0i1}, \eqref{Eq:geiIdgk} are restated as
 \begin{align}
 f^{i_0}_i\circ f^{i_0\mm i_1}_{i_0}&=f^{i_1}_i\circ f^{i_0\mm i_1}_{i_1};
 \label{Eq:fi0i1}\\
 \ge_i\circ\Id f^{i_k}_i&=\gf^{i_k}_i\ge_{i_k},\qq\text{for }k\in\set{1,2}.
 \label{Eq:geiIdfik}
 \end{align}
At this point, we have defined $f^j_i$, if $i\in I_{n+1}\setminus I_n$
and $j$ is an immediate predecessor of $i$. If $i\in I_{n+1}\setminus I_n$
and $j<i$, then the only possibility is to put
$f^j_i=f^{i_\gn}_i\circ f^j_{i_\gn}$, where $\gn<2$ is such that
$j\leq i_\gn$. For this to be possible, we need to verify that if
$j\leq i_0\mm i_1$, then $f^{i_0}_i\circ f^j_{i_0}=f^{i_1}_i\circ f^j_{i_1}$.
This follows from \eqref{Eq:fi0i1}, along with the following sequence of
equalities:
 \begin{align*}
 f^{i_0}_i\circ f^j_{i_0}&=
 f^{i_0}_i\circ f^{i_0\mm i_1}_{i_0}\circ f^j_{i_0\mm i_1}\\
 &=f^{i_1}_i\circ f^{i_0\mm i_1}_{i_1}\circ f^j_{i_0\mm i_1}\\
 &=f^{i_1}_i\circ f^j_{i_1}. 
 \end{align*}
At this point, we have defined $f^j_i$, if $i\in I_{n+1}\setminus I_n$
and $j<i$. We extend this definition by putting $f^i_i=\id_{R_i}$. The
verification of conditions (a)--(d) above is then straightforward.

Let $R$ be the direct limit of the $R_i$, with the transition maps $f^i_j$
for $i\leq j$ in $I$. Since the $\Idc$ functor preserves direct limits,
$\Idc R$ is isomorphic to the direct limit of the $\Idc R_i$, with the
transition maps $\Idc f^i_j$, for $i\leq j$ in $I$. By (c) and (d), it
follows that $\Idc R$ is isomorphic to the direct limit of all the $S_i$,
with the transition maps $\gf^i_j$, for $i\leq j$ in $I$; whence
$\Idc R\iso S$. This settles part (ii) of the statement of
Theorem~\ref{T:Repr}.

Let now $S$ be a distributive \jz-semilattice of cardinality at most
$\aleph_1$. Adjoin a new largest element to $S$, forming $T=S\uu\set{1}$.
Then $T$ is a distributive \jz-semilattice with a largest element and of
cardinality at most $\aleph_1$, thus, by what we just proved, it is
isomorphic to $\Idc R$, for some regular unital $K$-algebra $R$. Let
$\ge\colon\Idc R\to T$ be an isomorphism. We define an ideal $I$ of $R$, by
 \[
 I=\setm{x\in R}{\ge(RxR)\in S}.
 \]
Then $S$ is the image of $\Idc I$ under $\ge$. In particular, $\Idc I\iso S$.
Since $R$ is regular and $I$ is an ideal of $R$, $I$ is regular, see
\cite[Lemma 1.3]{GvnRR}.
\end{proof}

\begin{corollary}\label{C:Repr}
Let $D$ be an algebraic distributive lattice with at most $\aleph_1$ compact
elements. Then there exists a sectionally complemented modular lattice $L$
such that $\Con L\iso D$. Furthermore, if the largest element of $D$ is
compact, then one can take $L$ with a largest element.
\end{corollary}

\begin{proof}
By Theorem~\ref{T:Repr}, there exists a regular ring $R$ such that
$\Id R\iso D$, and $R$ is unital if the largest element of $D$ is compact.
Let $L$ be the lattice of all principal right ideals of $R$. Then $L$ is a
sectionally complemented modular lattice, see Section~\ref{S:basic}.
Note that if $D$ has a largest element, then $R$ is unital, thus $L$ has a
largest element. By \cite[Theorem~4.3]{Weurp}, $\Con L$ is isomorphic to
$\Id R$. Hence $\Con L\iso D$.
\end{proof}

\section{Problems and comments}\label{S:Pbs}

The main amalgamation result of this paper, Theorem~\ref{T:Amalg}, is
easily seen to imply that for every finite diagram $\E{D}$ of finite
Boolean \jz-semilattices, if $\E{D}$ is indexed by the \emph{square}
$\two^2$, then $\E{D}$ has a lifting, with respect to the $\Idc$ functor,
by regular rings. The analogue of this result in case $\E{D}$
is indexed by the \emph{cube}, $\two^3$, does \emph{not} hold, by the results
of \cite{TuWe}. A $\two^3$-indexed diagram $\E{D}$ is produced there,
that cannot be lifted, with respect to the $\Conc$ functor, by lattices with
permutable congruences. Since for a regular ring $R$, the ideal lattice of
$R$ is isomorphic to the congruence lattice of the principal right ideal
lattice $L$ of $R$, and since $L$ has permutable congruences (because it is
sectionally complemented), $\E{D}$ cannot be lifted, with respect to the
$\Idc$ functor, by regular rings as well.\smallskip

The regular rings obtained in the underlying construction of Bergman's
Theorem are \emph{locally matricial}, that is, direct limits of matricial
rings (a ring is \emph{matricial} over a division ring $K$, if it is
isomorphic to a direct product of finitely many matrix rings over $K$).
Quite to the contrary, the rings that we obtain in the proof of
Theorem~\ref{T:Repr} are not locally matricial. In fact, if $R$ is one of
those rings, then the monoid $V(R)$ of isomorphism classes of finitely
generated projective right $R$-modules is a \emph{semilattice}, as opposed
to the case of a locally matricial ring $R$, for which $V(R)$ is
\emph{cancellative}. Our proof cannot be extended to provide $R$ locally
matricial, because the finite-dimensional analogue of Theorem~\ref{T:Amalg}
fails. This suggests the following problem:

\begin{problem}\label{Pb:LocMat}
Let $D$ be a distributive algebraic lattice with at most $\aleph_1$ compact
elements. Does there exists a locally matricial ring $R$ such that
$\Id R\iso D$?
\end{problem}

Problem~\ref{Pb:LocMat} is equivalent to the particular instance of
\cite[Problem~1]{GoWe1} obtained by taking $|S|\leq\aleph_1$.

Because of Example~\ref{Eq:NonFin}, Theorem~\ref{T:Amalg} does not extend to
matricial algebras over a division ring. However, this does not rule out
\emph{a priori} the possibility of an extension of Theorem~\ref{T:Amalg} to a
sufficiently large subcategory of the category of matricial $K$-algebras and
unital homomorphisms of algebras.

\begin{problem}\label{Pb:LocFin}
Let $D$ be an algebraic distributive lattice with at most $\aleph_1$ compact
elements. Does there exist a \emph{locally finite}, complemented, modular
lattice $L$ such that $\Con L\iso D$?
\end{problem}

(A lattice $L$ is \emph{locally finite}, if every finitely generated
sublattice of $L$ is finite.)
It can be easily proved, by using some results of \cite{GoHa} and
\cite{WDim}, that a positive answer to Problem~\ref{Pb:LocFin} would imply a
positive answer to Problem~\ref{Pb:LocMat}. Conversely, if one obtains a
positive answer to Problem~\ref{Pb:LocMat} with $R$ locally matricial over a
\emph{finite} field $K$, then the lattice of principal right ideals of $R$ is
locally finite, thus providing a positive answer to Problem~\ref{Pb:LocFin}.

\end{document}